\documentclass[12pt]{extarticle}
\usepackage{amsmath, amsthm, amssymb, hyperref, color}
\usepackage[shortlabels]{enumitem}
\usepackage{graphicx}
\usepackage[all]{xypic}
\usepackage{makecell}
\usepackage[final]{pdfpages}
\setboolean{@twoside}{false}
\usepackage{pdfpages}
\usepackage{caption}
\usepackage{subcaption}
\usepackage{scalefnt}
\usepackage{verbatim}
\oddsidemargin \evensidemargin
\newtheorem{theorem}{Theorem}
\numberwithin{theorem}{section}
\newtheorem{proposition}[theorem]{Proposition}
\newtheorem{lemma}[theorem]{Lemma}

\newtheorem{definition}[theorem]{Definition}

\newtheorem{conjecture}[theorem]{Conjecture}
\newtheorem{algorithm}[theorem]{Algorithm}
\theoremstyle{remark}
\newtheorem{remark}[theorem]{Remark}
\newtheorem{example}[theorem]{Example}

\newcommand{\RR}{\mathbb{R}}
\newcommand{\QQ}{\mathbb{Q}}
\newcommand{\PP}{\mathbb{P}}
\newcommand{\CC}{\mathbb{C}}

\newcommand{\AAA}{\mathbb{A}}

 \date{}

\DeclareMathOperator{\rk}{rk} 
\DeclareMathOperator{\Spec}{Spec} 
\title{\textbf{Spaces of Sums of Powers and Real Rank Boundaries}}

\author{Mateusz Micha{\l}ek, Hyunsuk Moon}

\begin{document}

\maketitle
\thanks{Micha{\l}ek was supported by IP grant 0301/IP
3/2015/73 of the Polish Ministry of Science.}
\begin{abstract} \noindent
We investigate properties of Waring decompositions of real homogeneous forms. We study the moduli of real decompositions, so-called {\bf S}pace of {\bf S}ums of {\bf P}owers, naturally included in the {\bf V}ariety of {\bf S}ums of {\bf P}owers. Explicit results are obtained for quaternary quadrics, relating the algebraic boundary of ${\rm SSP}$ to various loci in the Hilbert scheme of four points in $\PP^3$. Further, we study the locus of general real forms whose real rank coincides with the complex rank. In case of quaternary quadrics the boundary of this locus is a degree forty hypersurface $J(\sigma_3(v_3(\mathbb{P}^3)),\tau(v_3(\mathbb{P}^3)))$.
\end{abstract}

\section{Introduction}
Let $f\in K[x_1,\dots,x_n]_d=S^d(V^*)$ be a homogeneous form of degree $d$ in $n$ variables over a field $K$. We study \emph{Waring decompositions}:
$$f=\sum_{i=1}^r\lambda_i (\sum_{j=1}^n a_{i,j}x_i)^d,\text{ where }\lambda_i,a_{i,j}\in K.$$
The smallest possible $r=\rk_K(f)$ is called the \emph{Waring rank} of $f$. 

Waring rank and decompositions have attracted attention of many mathematicians. On the one hand, their study is motivated by applications e.g.~in computer sciences. On the other hand these problems are related to beautiful mathematics: geometry (through secant varieties) \cite{Zak}, representation theory (through homogeneous varieties) \cite{LM}, algebra (through apolar ideals and resolutions) \cite{RS}, moduli spaces \cite{Muk1} and many more \cite{LandsbergTensor}.

The main questions that motivated us are:
\begin{enumerate}
\itemsep0em
\item How to find a Waring decomposition for a general form $f$?
\item What is the geometry (moduli) of all decompositions?
\item How does the answer to first two questions depend on $K=\RR$ or $K=\CC$?
\item When a real form admits a real Waring decomposition with $r=\rk_\CC(f)$?
\end{enumerate}

To investigate the geometry of the Waring decompositions one defines the {\bf V}ariety of {\bf S}ums of {\bf P}owers \cite{RS}:
$${\rm VSP}(f):=\overline{\{(\ell_1,\dots,\ell_r):f=\sum\lambda_i\ell_i^d, [\ell_i]\in \PP(V^*)\}}.$$
Here, the closure is taken in the (smoothable component of the) Hilbert scheme of $r$ points in $\PP(V^*)$ and is usually considered for $K=\CC$. Seminal works of Mukai, Ranestad, Schreyer \cite{Muk1, Muk2, Muk3, RSpolar, RS} and others provide descriptions of ${\rm VSP}$'s in special cases. Following \cite{MMSV}, we investiagate a semialgebraic subset of the real locus of ${\rm VSP}$, {\bf S}pace of {\bf S}ums of {\bf P}owers, corresponding to those decompositions in which all $\ell_i$ are real. We note that a real point of ${\rm VSP}$ does \emph{not} have to correspond to fully real decomposition. The geometry and topology of the real locus of ${\rm VSP}$ is another very interesting subject, studied e.g.~in \cite{KS}. 

The first part of the paper is devoted to the algebraic boundary $\partial_{alg} {\rm SSP}(f)$ -- the Zariski closure of difference of the Euclidean closure of ${\rm SSP}$ minus the interior of the closure.

We start by defining the anti-polar form $\Omega(f)$ for a general form of even degree in Definition \ref{def:anti}. This is a generalization of the dual quadric. As we show, in several cases $\Omega(f)$ governs the nonreduced structure of apolar schemes -- cf.~Proposition \ref{prop:degeneration}.  
This allows us to provide an explicit description of $\partial_{alg}{\rm SSP}(f)$, when $f$ is a quaternary quadric \ref{th:Quat Quad}, extending the results obtained for ternary quadrics \cite[Section 2]{MMSV}. In the ternary case, due to results on Gorenstein resolutions, one could provide a description of $\partial_{alg}{\rm SSP}(f)$ in terms of hyperdeterminants. We did not find a straightforward generalization to quaternary case, however our methods may be applied to ternary case. They provide a simpler, more explicit (but less intrinsic) description \ref{sec:TerQuad}. Further, we study the geometry of $\partial_{alg}{\rm SSP}(f)$ in detail, over $\QQ$ and over $\CC$, relating it to the geometry of the Hilbert scheme. In particular, we discuss how the ${\rm VSP}$ intersects various loci of the Hilbert scheme. 

We note that our approach towards $\partial_{alg}{\rm SSP}(f)$ allows a description of intersection with subvarieties of ${\rm VSP}(f)$. This is useful, when the ${\rm VSP}(f)$ itself is a large variety. We apply it for quinary quadrics \ref{sec:QuinQuad}.

In the second part of the article we study the locus of real general forms for which the real rank equals the complex rank. In general, the complex rank is known by Alexander-Hrischowitz theorem \cite{AH}. It follows that there exists a Zariski dense semialgebraic set $\mathcal{R}_{n,d}$ of real forms with such real rank. The \emph{real rank boundary} $\partial_{alg}\mathcal{R}_{n,d}$ is defined as the Zariski closure of the topological boundary of $\mathcal{R}_{n,d}$. 

We present an implemented, fast, deterministic algorithm, that for a general quaternary cubic $f$ returns its unique Waring decomposition. This is an interesting example, one of the only two identifiable that is not binary \cite{GaluppiMella}, studied both classically \cite{Clebsch} and in modern context \cite{OttavianiOeding}. Our algorithm can be used not only in applications, but working parametrically over the field $K(t)$ it allows to provide a description of $\partial_{alg}\mathcal{R}_{4,3}$. We obtain Proposition \ref{prop:R43};

The join variety $J(\sigma_3(v_3(\mathbb{P}^3)),\tau(v_3(\mathbb{P}^3)))$ of the third secant variety of the third Veronese of $\PP^3$ and the tangential variety is an irreducible hypersurface of degree $40$ in $\PP^{19}$. It equals $\partial_{alg}(\mathcal{R}_{4,3})$.

The geometry of $\partial_{alg}\mathcal{R}_{n,d}$ is also related to the classical work of Hilbert on cones of sums of squares and nonnegative polynomials \cite{Hilbert}. This allows us to provide a description of one component of $\partial_{alg}\mathcal{R}_{4,4}$ as the dual of a variety of quartic symmetroids in Proposition \ref{prop:R44}.

In the Appendix \ref{sec:App} we present several algorithms used to:
\begin{itemize}
\itemsep0em
\item prove that a given variety is irreducible (or reducible) over $\CC$,
\item compute the real locus of a variety,
\item describe a shape of the resolution of a generic member of a family of ideals.
\end{itemize}
\section*{Acknowledgments}
We would like to thank Joachim Jelisiejew for interesting discussions about the geometry of the Hilbert scheme. We thank Grigoriy Blekherman and Rainer Sinn for pointing us towards the results on positive real rank. Micha{\l}ek was supported by the Foundation for Polish Science (FNP) and is a member of AGATES group. Part of the research was realized during research visits at FU Berlin and RIMS in Kyoto - we express our gratitude to Klaus Altmann, Takayuki Hibi and Hiraku Nakajima for their hospitality.
\section{Boundaries of Spaces of Sums of Powers}\label{sec:Boud of SSP}

\begin{definition}[The anti-polar $\Omega(f)$]\label{def:anti}{\rm
Consider a homogeneous polynomial $f\in S^{2d}(V^*)$ of degree $2d$. It induces, through the middle catalecticant, a linear map $A_f:S^d(V)\rightarrow S^d(V^*)$. Suppose that $A_f$ is an isomorphism (which holds for generic $f$). The inverse map $A_f^{-1}$ defines a polynomial $\Omega(f)\in S^{2d}(V)$ as follows. For $x\in V^*$ we define $\Omega(f)(x):=<x^d, A_f^{-1}(x^d)>$, where $<\cdot,\cdot>$ denotes the perfect paring of $S^d(V)$ and $S^d(V^*)=S^d(V)^*$.}
\end{definition}
\begin{remark}\label{rem:explicit}{\rm
Explicitly for a fixed base on $V$, to evaluate $\Omega(f)$ on $x$ one multiplies the inverse (or adjoint up to scalar) of the middle catalecticant from left and from right by the vector that evaluates all monomials of degree $d$ on $x$.}
\end{remark}
The following proposition explains the name of $\Omega(f)$, relating it to the anti-polar quartic defined in \cite[Theorem 4.1]{MMSV}.
\begin{proposition}
For a homogeneous polynomial $f\in S^{2d}(V^*)$ let us denote by $C(f)$ its middle catalecticant. If $C(f)$ is invertible, then (up to scalar) we have:
$$\Omega(f)(l) \,\, := \,\, \,{\rm det}\bigl(C(f + \ell^{2d}) \bigr)\, -\, {\rm det}\bigl(C(f) \bigr),$$
where $l\in V^*$.
\end{proposition}
\begin{proof}
As the Definition \ref{def:anti} and the one in the proposition are intrinsic, we may fix a basis and assume $l$ is a basis element. Then $C(l^{2d})$ is simply given by a matrix with one nonzero entry on the diagonal. Hence, ${\rm det}\bigl(C(f + l^{2d}) \bigr)- {\rm det}\bigl(C(f) \bigr)$ equals the complimentary minor to the nonzero entry. This exactly agrees with Definition \ref{def:anti} by Remark \ref{rem:explicit}.
\end{proof}
\begin{example}{\rm
If $f$ is a quadric of full rank, than $\Omega(f)$ is simply the dual quadric.}
\end{example} 
\begin{remark}
{\rm It is important to note that neither the inverse nor the adjoint of $A_f$ is the middle catalecticant of $\Omega(f)$ - cf.~\cite{Dol1} and \cite[Remark 1.4.1]{Dol2}, contrary to the case of quadrics. In fact, it is often the case that $A_f^{-1}$ is not a middle catalecticant of any equation of degree $2d$ - see e.g.~\cite[Proposition 7.1]{MSUZ}.}
\end{remark}
The following proposition is based on results from \cite{RSpolar}.
\begin{proposition}\label{prop:degeneration}
Let $f\in S^{2d}(V^*)$ be such that the middle catalecticant $A_f$ is nondegenerate. Suppose $S\subset \PP(V)$ of length equal to the rank of $A_f$ is apolar to $f$. Then $S$ has a nonreduced structure at a point $l\in S$ if and only if $\Omega(f)(l)=0$. 
\end{proposition}
\begin{proof}
Consider a new scheme $\tilde S$ defined by:
$$I_{\tilde S}:=I_S:l^{\perp},$$
where $l^\perp$ is the maximal ideal defining $l$. As $A_f$ is nondegenerate we know that $(I_S)_d=\emptyset$. On the other hand, by degree count we know that there exists $g\in (I_{\tilde S})_d$. We may determine the $g$ as follows:
$$g\in I_{\tilde S}\Leftrightarrow gl^\perp\subset I_S\Rightarrow g(l^\perp)_d\subset I_S\Rightarrow (g(l^\perp)_d)(f)=0\Leftrightarrow g((l^\perp)_d(f))=0.$$
But $(l^\perp)_d(f)=A_f(l^\perp_d)$, so
$
g\in I_{\tilde S}\Rightarrow g(A_f(l^\perp_d))=0.
$
As $A_f$ is nondegenerate and $l^\perp_d$ is a hypersurface in $S^dV^*$ this determines $g$ uniquely up to scalar and the above implication is an equivalence. We obtain that $S$ is reduced at $l$ if and only if $v_d(l)\not \in A_f(l^\perp_d)$ if and only if $A_f^{-1}(v_d(l))\not\in (l^\perp)_d$, where $v_d$ is the $d$-th Veronese (evaluating $g$ as a polynomial on $l$ is the same as evaluating $g$ as a linear map on $v_d(l)$). Hence, $S$ is nonreduced at $l$ if and only if $<v_d(l),A_f^{-1}(v_d(l))>= 0$. 
\end{proof}
The following definitions will be useful in the study of ${\rm SSP}$.
\begin{definition}[${\rm VNSP}$, $\mathfrak F$]
For a given form $f\in S^d(V^*)$ we define the {\bf V}ariety of {\bf N}onreduced {\bf S}ums of {\bf P}owers as a subscheme with the reduced structure of ${\rm VSP}(f)$ \mbox{corresponding} to the locus of nonreduced schemes. In other words:
$${\rm VNSP}(f)={\rm VSP}(f)\setminus \{(\ell_1,\dots,\ell_r):f=\sum\lambda_i\ell_i^d, [\ell_i]\in \PP(V^*)\},$$
or equivalently it parametrizes smoothable, nonsmooth apolar schemes. 

As a moduli space ${\rm VSP}$ comes with a universal family $\pi:{\rm VSP}(f)\times \PP(V^*)\supset \mathfrak{F}\rightarrow {\rm VSP}(f)$, where the fiber over a given point of ${\rm VSP}(f)$ equals the corresponding apolar scheme.
\end{definition}

Let us describe how Proposition \ref{prop:degeneration} may be used to find the boundry of ${\rm SSP}$ inside the ${\rm VSP}$ in special cases. Let $f\in S^{2d}(V^*)$ have a nondegenerate middle catalecticant. Further assume that the rank of the catalecticant equals the generic rank in $S^{2d}(V^*)$.
\begin{enumerate}
\itemsep0em
\item Assume we are given the universal family $\pi:{\rm VSP}(f)\times \PP(V^*)\supset \mathfrak{F}\rightarrow {\rm VSP}(f)$.
\item Let $\mathfrak{B}={\rm VSP}(f)\times V(\Omega(f))\subset {\rm VSP}(f)\times \PP(V^*)$.
\item The algebraic boundary of  ${\rm SSP}(f)$ inside ${\rm VSP}(f)$ is contained in $\pi(\mathfrak{B}\cap \mathfrak{F})$. 
\end{enumerate}
We note that set-theoretically $\pi(\mathfrak{B}\cap\mathfrak{F})$ coincides with ${\rm VNSP}(f)$.
The construction above follows from the fact that a real decomposition can change to a complex one only by passing trough a nonreduced scheme. 
As a consequence we obtain the following lemma that will allow us to find defining equations of $\partial_{alg}{\rm SSP}(f)$.

\begin{lemma}\label{lem:projection and boundary of SSP}
Suppose that ${\rm SSP}(f)\neq\emptyset$. If the top dimensional component of $\pi(\mathfrak{B}\cap\mathfrak{F})$ is irreducible, then (with reduced structure) it coincides with $\partial_{alg}{\rm SSP}(f)$.
\end{lemma}

In principal, this method could provide the description for the boundary of the ${\rm SSP}$ inside the ${\rm VSP}$ for quaternary quartics. Unfortunately, in this case the ${\rm VSP}$, which is a $5$-fold eludes an explicit description. Citing Ranestad it is "one of the most interesting outstanding problems on ${\rm VSP}$'s."



\subsection{Quadrics}
We now apply Proposition \ref{prop:degeneration} to explicitly obtain the boundary of ${\rm SSP}$ for quadrics $f$ in up to $n\leq 5$ variables, i.e.~in all cases when ${\rm VSP}$ is smooth. For $n=2,3$ this was achieved in \cite{MMSV}. It that case, as the codimension of the apolar ideal $f^\perp$ was at most three, one could apply the classical results of Buchsbaum-Eisenbud on resolutions of Gorenstein schemes \cite{BuchsbaumEisenbud} and define the boundary by an appropriate hyperdeterminant. In case $n>3$ we could still take the resolution of $f^\perp$, however explicit results using this technique seem much harder. Instead we follow Proposition \ref{prop:degeneration}. 
\begin{lemma}\label{lem:signaturethe sameboundary}
Fix two real forms $f_1,f_2$ with (both) ${\rm SSP}$'s nonempty. Suppose $f_1$ can be obtained from $f_2$ by a \emph{complex} change of coordinates. If the top dimensional component of ${\rm VNSP}(f_1)\subset {\rm VSP}(f_1)$ is irreducible, then $\partial_{alg} {\rm SSP}(f_1)$ is isomorphic (as a complex algebraic variety) to $\partial_{alg} {\rm SSP}(f_2)$.
\end{lemma}
\begin{proof}
The isomorphism between $f_1$ and $f_2$ provides an isomorphism between their ${\rm VSP}$'s. We show that this isomorphism is also an isomorphisms of the algebraic boundaries.
By the assumption the boundary is nonempty in both cases, hence of codimension one. However, it has to be contained in ${\rm VNSP}$, which is irreducible and preserved by the isomorphism.
\end{proof}
\begin{remark}
Note that we need the assumption that ${\rm SSP}$'s are nonempty. We have ${\rm SSP}(x^2+y^2)={\rm VSP}(x^2+y^2)$, hence the algebraic boundary is empty. On the other hand, a quadric with a different signature (but of course isomorphic over $\CC$) satisfies
${\rm SSP}(x^2-y^2)\subsetneq{\rm VSP}(x^2-y^2)$ and the algebraic boundary consists of two points.
\end{remark}
We will apply Lemma \ref{lem:signaturethe sameboundary} to quadrics of different signature. Our aim is to describe the algebraic boundary of the ${\rm SSP}$.
\subsubsection{Ternary Quadrics}\label{sec:TerQuad}
The case of ternary quadrics is well-understood \cite[Section 2]{MMSV}. We use it as a warm-up. Fix $f=x_1x_3+x_2^2$, which up to real isomorphism is the only indefinite quadric. The ${\rm VSP}(f)$ is a smooth Fano 3-fold $V_5$ -- quintic del Pezzo threefold -- admitting a realization as an intersection $G(3,5)\cap \PP^6$. The boundary $\partial_{alg}{\rm SSP}(f)$ is given by a special hyperdeterminant, which turns out to be a degree $20$ polynomial in $6$ variables with $13956$ terms. 

We show how this polynomial simplifies, if we work in a local affine patch of ${\rm VSP}(f)$ described in \cite{RSpolar}. Indeed, computing $\pi(\mathfrak{B}\cap\mathfrak{F})$ in Appendix \ref{sec:AppTerQuad} we obtain a quartic surface:
$$27a^2-32b^3+36abc-4b^2c^2+4ac^3.$$
Its singular locus is a smooth curve -- a complete intersection of a quadric and cubic surface. It corresponds to the locus where the apolar scheme is local, i.e.~all three points coming together. As the quartic is irreducible it follows that it coincides with $\partial_{alg}{\rm SSP}(f)$.
\begin{figure}
\hskip-0.3cm
\includegraphics[scale=0.33]{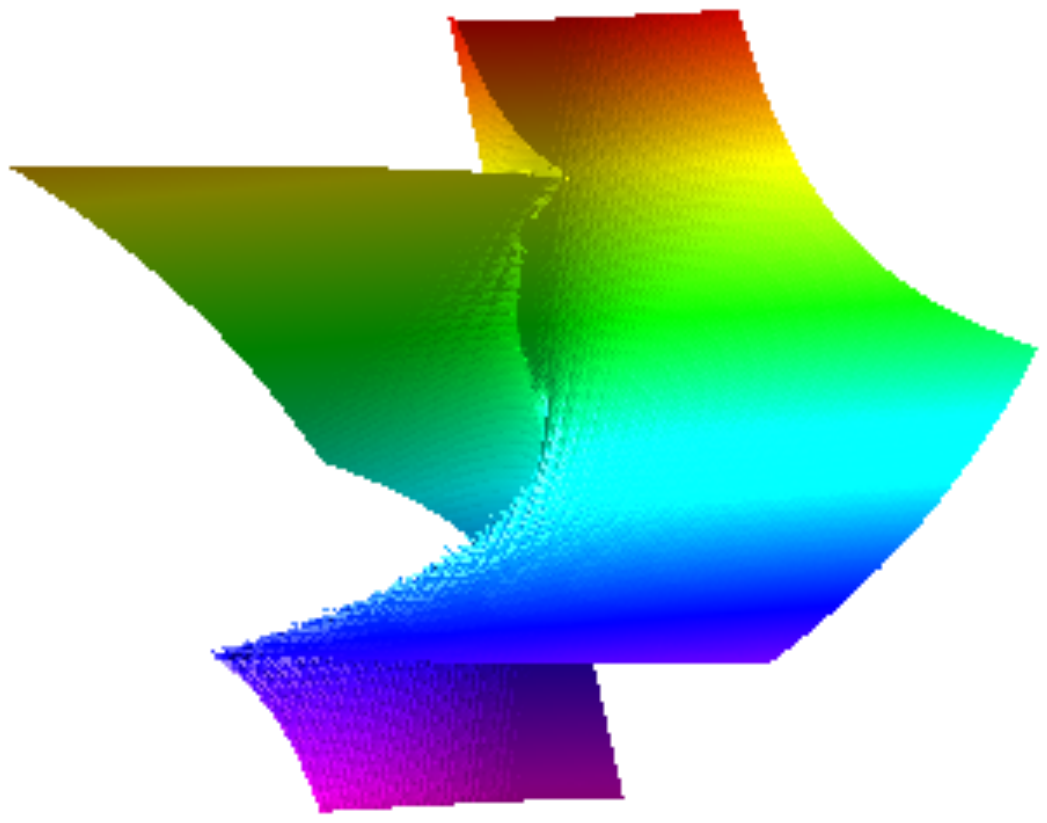}\hskip-2.5cm
\includegraphics[scale=0.33]{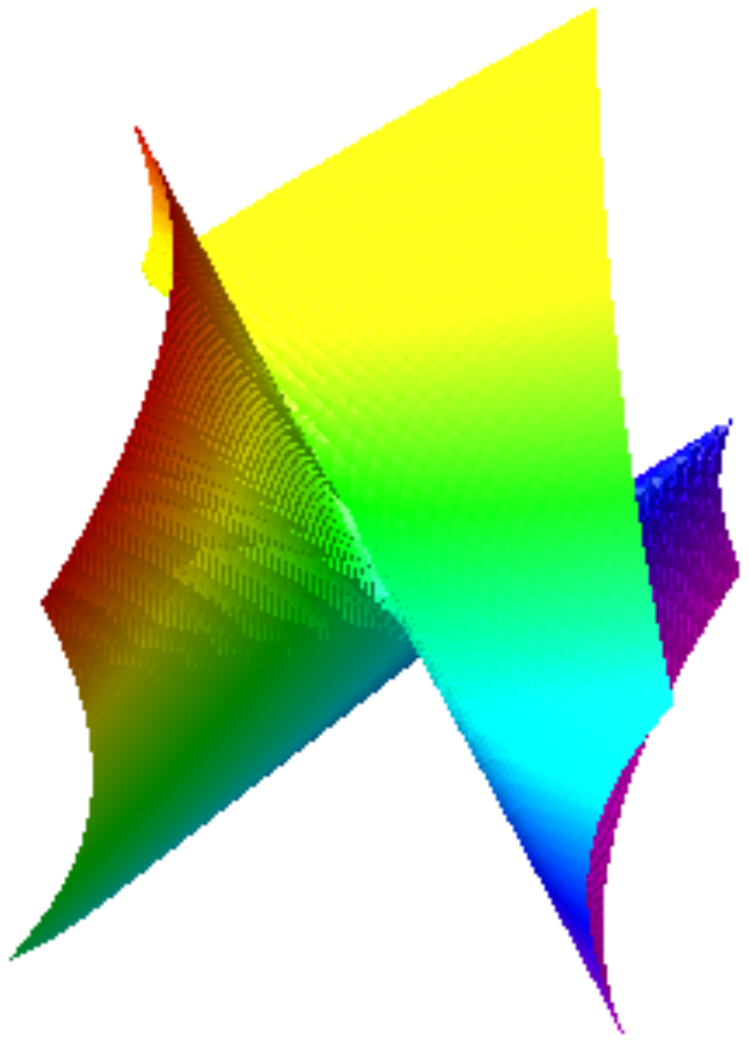}\hskip-2.5cm
\includegraphics[scale=0.33]{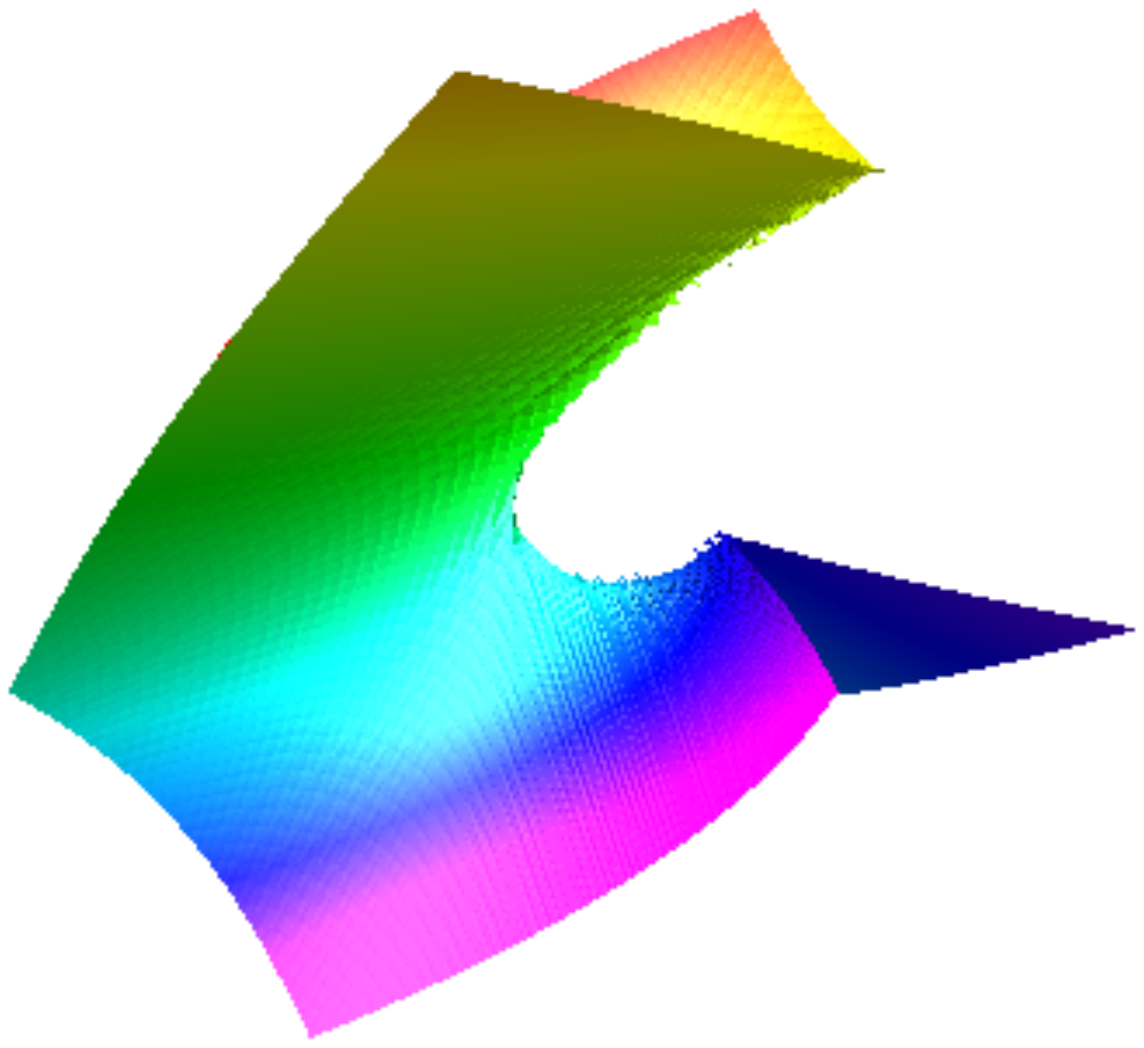}
\vskip-3cm
Three pictures of the affine patch of $\partial_{alg}{\rm SSP}(f=x_1x_3+x_2^2)$. The three-dimensional ambient affine space represents the quintic Fano threefold - moduli of schemes apolar to $f$ of length three. The surface corresponds to those schemes that are supported at two points: one smooth, one with reduced structure $\Spec\CC[x]/(x^2)$. The curve in the singular locus represents apolar schemes isomorphic to $\Spec\CC[x]/(x^3)$. The surface divides the threefold into two regions: one is the ${\rm SSP}(f)$ corresponding to fully real decomposition, the other corresponds to decompositions where one linear form is real and the other two are conjugate.
\end{figure}
\subsubsection{Quaternary Quadrics}
Let $f=x_1x_4+x_2^2+x_3^2$. The following theorem computes the boundary of the ${\rm SSP}$ relating its geometry to the Hilbert scheme and the geometry of the apolar schemes. Ranestad and Schreyer \cite{RS} provide an explicit local description of the variety ${\rm VSP}(f)$. We will be working locally on such affine patches.
\begin{theorem}\label{th:Quat Quad}\
\begin{enumerate}
\itemsep0em
\item The variety $\mathfrak{Y}=\partial_{alg}{\rm SSP}(f)\subset {\rm VSP}(f)$ is irreducible over $\CC$, of dimension $5$ and equals $\pi(\mathfrak{B}\cap\mathfrak{F})$.
\item The singular locus of $\mathfrak{Y}$ has two four dimensional components over $\QQ$: $\mathfrak{Y}_{3,1}$ and $\mathfrak{Y}_{2,2}$. The general point of $\mathfrak{Y}_{3,1}$ (resp.~$\mathfrak{Y}_{2,2}$) corresponds to a nonreduced scheme with two support points; one of which has multiplicity $3$ (resp.~$2$), the other has multiplicity $1$ (resp.~$2$).
\item Over $\CC$, $\mathfrak{Y}_{2,2}$ however has two irreducible components; $\mathfrak{Y}_{2,2}'$ and $\mathfrak{Y}_{2,2}''$ intersecting along surface (identified later with $\mathfrak{Y}_{4,sing}$).
\item The two components $\mathfrak{Y}_{3,1}$, $\mathfrak{Y}_{2,2}$ intersect along a three dimensional threefold $\mathfrak{Y}_{4}$, irreducible over $\QQ$. It is the singular locus of $\mathfrak{Y}_{3,1}$.
\item Over $\CC$, $\mathfrak{Y}_{4}$ has two components, corresponding to intersections of $\mathfrak{Y}_{2,2}'$ and $\mathfrak{Y}_{2,2}''$ with $\mathfrak{Y}_{3,1}$.
\item The general point of $\mathfrak{Y}_{4}$ corresponds to a nonreduced, local scheme of length four, isomorphic to ${\rm Spec} \CC[x]/(x^4)$. 
\item The singular locus $\mathfrak{Y}_{4,sing}$ of $\mathfrak{Y}_{4}$ coincides with the singular locus of $\mathfrak{Y}_{2,2}$. It is a two dimensional smooth surface. The general point of $\mathfrak{Y}_{4,sing}$ corresponds to a nonreduced, local scheme isomorphic to $$\Spec \CC[x,y]/(x^2-y^2,xy)$$.
\item The locus of real points of $\mathfrak{Y}_{4}$ and $\mathfrak{Y}_{4,sing}$ coincide.
\end{enumerate}
\end{theorem}
By Lemma \ref{lem:signaturethe sameboundary} the same result holds for all indefinite, full rank quadrics. 
\begin{proof}
By \cite[Theorem 1.1]{RSpolar} ${\rm VSP}(f)$ is a smooth, six dimensional variety. Ranestad and Schreyer provide an explicit local description of this variety and the universal family $\mathfrak{F}$ in a Macaulay 2 package VarietyOfPolarSimplices.m2 \cite{Mac2:varPS}. We extend their defining equations by $\Omega(f)$, obtaining $\mathfrak{B}\cap\mathfrak{F}\subset {\rm VSP}(f)\times \PP^3$. Eliminating the variables corresponding to $\PP^3$ automatically is not possible (due to a large ambient dimension of the ${\rm VSP}$). However, one may find automorphisms of the (local) description of the ${\rm VSP}$ that reduce the number of variables. 
Afterwards we perform elimination, obtaining explicitly the defining equation of $\pi(\mathfrak{B}\cap\mathfrak{F})$ (defined over $\QQ$)-see Appendix \ref{app:Quat Quad}. We check that it defines a prime ideal, over $\CC$, as follows -- details of implementation are presented in Appendix \ref{App:Boundary prime}. We fix four linear forms defined over $\QQ$ and intersect them with the $\pi(\mathfrak{B}\cap\mathfrak{F})$ obtaining a curve $C$. It is enough to show that $C$ is irreducible over $\CC$, as $\pi(\mathfrak{B}\cap\mathfrak{F})$ was a hypersurface -- in particular equidimensional. We project, not changing the degree, until $C$ becomes a plane curve, that is defined by $g_C$. We prove it is irreducible as follows.
We consider all possible factorizations of $g_C$ into a product of degrees $4+4$, $3+5$, $2+6$, $1+7$ with coefficients given by variables. Comparing all coefficients we obtain ideals, that equal the whole ring.
By the previous discussion this proves the statement $1$. 

The computation of the singular locus of $\mathfrak{Y}$ and its decomposition over $\QQ$ is now straightforward. The two components $\mathfrak{Y}_{3,1}$ and $\mathfrak{Y}_{2,2}$ are distinguished by the dimension of their singular locus, which is respectively $3$ and $2$. There are several ways to prove that $\mathfrak{Y}_{3,1}$ corresponds to schemes of type $(3,1)$ and $\mathfrak{Y}_{2,2}$ of type $(2,2)$. One can restrict the family $\mathfrak{F}$ to $\mathfrak{Y}_{3,1}$ and intersect it with $\Omega(f)$. By Proposition \ref{prop:degeneration} in case of schemes of type $(3,1)$ these yields a family of local schemes (as only one of the points was nonreduced), while the schemes of type $(2,2)$ provide two distinct points. Statements $2$ and $4$ follow.

To prove the statement $3$ we project $\mathfrak{Y}_{2,2}$ obtaining a hypersurface $H$ of degree $4$. If $\mathfrak{Y}_{2,2}$ were irreducible, then $H$ would have to be irreducible. However, the equation defining $H$ decomposes as a product of two quadrics. 

It may seem not clear why $\mathfrak{Y}_{4}$ - corresponding to all four points coming together - allows further degeneration. The reason is the geometry of the punctual Hilbert scheme of schemes of length four. The smoothable component is irreducible and consists of alignable schemes with a general point corresponding to the aligned scheme $\Spec \CC[x]/(x^4)$. However, not all schemes are aligned. The scheme $\Spec\CC[x,y]/(x^2-y^2,xy)$ is local, Gorenstein, but has a two dimensional tangent space. A very explicit degeneration of four points to this scheme is presented e.g.~in \cite[p.~11]{LanMich}. The main difference between $\CC[x]/(x^4)$ and $\CC[x,y]/(x^2-y^2,xy)$ is that in the first case the four points degenerate on a line, while in the second they come from distinct, linearly independent directions.

More explicitly, for a scheme corresponding to a point of $\mathfrak{Y}_{4}$ one can find an explicit projection to a line that preserves its degree - a feature possible only in case of $\Spec\CC[x]/(x^4)$. Further, one can find explicitly points of $\mathfrak{Y}_{4,sing}$ and directly prove that the associated schemes are isomorphic to $\Spec\CC[x,y]/(x^2-y^2,xy)$.   

The last statement follows by computation presented in Appendix \ref{App:RealPoints}. Precisely, we show that one equation in the ideal of $\mathfrak{Y}_{4}$ is of the form $f_1^2+f_2^2$ and $I(\mathfrak{Y}_{4})+(f_1,f_2)=I(\mathfrak{Y}_{4,sing})$.
\end{proof}
\begin{remark}
In the given embedding, $\mathfrak{Y}_{2,2}$ is of degree $4$. Thus, $\mathfrak{Y}_{2,2}'$ and $\mathfrak{Y}_{2,2}''$ are of degree $2$ each. As they are also of codimension 2, each of them should be defined by a quadric in a hyperplane section. Using arithmetics in a finite field we have checked that they are smooth.
\end{remark}
\begin{remark}
As the locus in the Hilbert scheme corresponding to schemes isomorphic to $\CC[x]/(x^4)$ is irreducible, it is surprising that $\mathfrak{Y}_{4}$ is reducible over $\CC$. It would be interesting to know if it could happen for generic form $f$ in $n+1$ variables of rank $r$ that:
\begin{itemize}
\itemsep0em
\item the codimension one subvariety of the smoothable component of the Hilbert scheme of $r$ points in $\PP^n$ corresponding to nonreduced, Gorenstein schemes is irreducible, 
\item the boundary $\partial_{alg}{\rm SSP}(f)$ is reducible and positive dimensional.
\end{itemize}
There is always a divisor in the smoothable component of the Hilbert scheme, which is corresponding to $r-2$ reduced points and a multiplicity $2$ point. In some cases, however, there are other \emph{divisors} corresponding to nonreduced Gorenstein schemes. For more information we refer the reader to \cite[Theorem 1]{Iarrobino} and further consequences of this fact are discussed in \cite[Example A.17]{BJJM}, \cite[p.~7]{LanMich}. However, even in such a case we do not know if $\partial_{alg}{\rm SSP}(f)$ could be reducible (although we expect it).
\end{remark}

\begin{remark}
All apolar schemes  that we obtain must be Gorenstein by \cite[Lemma 2.3]{BB14}. Indeed, otherwise there would exist a Gorenstein scheme of smaller length that would be apolar to the form. In this range (e.g.~in ambient dimension at most three or if length is at most ten) all Gorenstein schemes are smoothable.
Hence, the form would have smaller border rank, i.e.~would not be generic.
\end{remark}

\subsubsection{Quinary Quadrics}\label{sec:QuinQuad}
In this case the ${\rm VSP}$ is a ten dimensional smooth variety. As before we may work locally in an affine patch described by Ranestad and Schreyer. However, in this case explicit computation of $\pi(\mathfrak{B}\cap\mathfrak{F})$ is much harder, as the elimination of variables is both time and memory consuming. 
Instead we apply Lemma \ref{lem:projection and boundary of SSP}. We proceed as follows:
\begin{enumerate}
\itemsep0em
\item Eliminate by hand possibly many variables, using explicit isomorphisms of the ambient affine space of the affine patch of $\mathfrak{B}\cap\mathfrak{F}$.
\item Consider the compactification $\mathfrak{C}$ of the affine patch of $\mathfrak{B}\cap\mathfrak{F}$ in a projective space.
\item Fix a $\PP^2$ in a projective space that is a compactification of the affine patch of ${\rm VSP}$.
\item Restrict the family $\mathfrak{C}$ to the given $\PP^2$ and project obtaining $\mathfrak{C'}$. 
\item Check that $\mathfrak{C'}$ is a reduced, irreducible curve of degree $10$. 
\end{enumerate}
As $\mathfrak{C'}$ is reduced, irreducible and of degree $10$, it follows that the maximal dimensional component of $\pi(\mathfrak{B}\cap\mathfrak{F})$ must be irreducible and defined by an unique polynomial of degree $10$. A simplified, explicit computation is presented in Appendix \ref{App:QuinaryQuadrics}. The code uses finite fields to make the computations fast (performing formal computation takes several hours, but is possible with the same code after changing the field). 
By Lemma \ref{lem:projection and boundary of SSP} we obtain the following proposition.
\begin{proposition}
In the affine space $\AAA^{10}\subset{\rm VSP}(f=x_1x_5+x_2^2+x_3^2+x_4^2)$ the algebraic boundary $\partial_{alg}{\rm SSP}(f)$ is an irreducible hypersurface of degree $10$.
\end{proposition}
We further note that the curve $\mathfrak{C'}$ representing $\partial_{alg}{\rm SSP}(f)$ has $30$ singular points.

\section{The Real Rank Boundary}

\subsection{Quaternary Cubics}
A general quarternary cubic $f\in\RR[x_1,x_2,x_3,x_4]$ has complex rank $R(4,3)=5$ and the decomposition is uniquely given by
\begin{equation}
\label{eq:ffive}
f=\ell_1^3+\ell_2^3+\ell_3^3+\ell_4^3+\ell_5^3
\end{equation}
which was claimed by Sylvester in 1851, and proved by Clebsch in 1861 \cite{Clebsch}. It is known as Sylvester Pentahedral Theorem.
Similarly to \cite[Section 5]{MMSV}, we propose the following algorithm to compute the five linear forms $\ell_i$.

\begin{algorithm}\label{Decomposition of cubic}
\underbar{\rm Input}: A general quarternary cubic $f$.
\underbar{\rm Output}: The decomposition {\rm (\ref{eq:ffive})}.

\begin{enumerate}
\itemsep0em
\item Compute the apolar ideal $f^{\perp}$.
It is generated by six quadrics $g_1,g_2,g_3,g_4,g_5,g_6$.
\item Compute the syzygies of $f^{\perp}$. 
Find the five linear syzygies $l_{ij},i=1,\ldots,5,j=1,\ldots,6$ on the quadrics satisfying
$$\sum_{j=1}^6 l_{ij}g_j=0\rm{\quad for\quad all} \quad i=1,\ldots,5 $$
\item  Compute a vector $(c_1,c_2,c_3,c_4,c_5,c_6) \in \RR^6 \backslash \{0\}$
that satisfies $$c_1 l_{i1} + c_2 l_{i2} + c_3 l_{i3} + c_4 l_{i4} + c_5 l_{i5} + c_6 l_{i6} = 0$$ for all $i=1\ldots 5$.
\item Let $J$ be the ideal  generated by the quadrics  $\,c_6 g_1 - c_1 g_6$, 
$\,c_6 g_2 - c_2 g_6\,$,$\ldots$, $\,c_6 g_5 - c_5 g_6$. Compute the variety $V(J)$ in $\PP^3$.
It consists precisely of the points dual to $\ell_1,\ell_2,\ldots,\ell_5$.
\end{enumerate}
\end{algorithm}


To prove the correctness of this algorithm, we need the following lemma.

\begin{lemma}\label{correctness}
For a general quarternary cubic $f\in \RR[x_1,x_2,x_3,x_4]$, the minimal free resolution of its apolar ideal is given by
\begin{equation}\label{resolution of apolar}0\longrightarrow S\longrightarrow S^6 \stackrel{M}{\longrightarrow} S^{10} \stackrel{N}{\longrightarrow} S^6 \longrightarrow S \longrightarrow 0.
\end{equation}
Further, the map $M$ can be represented by $6\times10$ matrix with two blocks
\begin{equation}\label{matrix representation}M=[\quad L\quad|\quad Q\quad]\end{equation}
where $L$ is a $6\times 5$ matrix with linear forms and $Q$ is a $6\times5$ matrix with quadric forms.
Moreover, there exists a row operation that makes the last row of $L$ equal to $0$. It distinguishes a unique five dimensional subspace of quadrics in $f^\perp$ with five linear syzygies.
\end{lemma}
\begin{proof}
Since the apolar ideal of a quarternary cubic is Gorenstein, codimension four with six quadric generators, as in \cite{MR}, the resolution has the form (\ref{resolution of apolar}). So we only need to show that the matrix representation $M$ has the form (\ref{matrix representation}).

To prove this, consider a family of cubics in $A[x_1,x_2,x_3,x_4]$ where $A=\RR[a_0,\ldots,a_3]$;
$$f=a_0x_1^3+a_1x_1^2x_2+a_2x_1^2x_3+a_3x_1^2x_4+11x_1x_2^2-12x_1x_2x_3+7x_1x_2x_4+32x_1x_3^2-28x_1x_3x_4+11x_1x_4^2$$
$$+8x_2^3-13x_2^2x_3+34x_2^2x_4+19x_2x_3^2-38x_2x_3x_4+16x_2x_4^2+7x_3^3-41x_3^2x_4+7x_3x_4^2+13x_4^3$$
where sixteen of its coefficients are fixed.
After applying $GL_{\RR}(4)$-action on this family, it is dominant on the space of cubics, which we prove using an easy computation of a Jacobian. Hence, it is enough to prove the lemma holds for this family.

Using the code in Appendix \ref{sec:AppQuartCub}, we can show that the apolar ideal $f^\perp$ has six quadric generators  $g_1,\ldots,g_6$ in $S$. Also, there are exactly five independent linear syzygies  $l_{ij},i=1\ldots 5$, and a unique nonzero vector $(c_1,\ldots,c_6)\in A^6$, satisfying  
$$\sum_{j=1}^6 c_j l_{ij} = 0, \quad i=1,\ldots,5.$$

Hence, the first syzygy matrix $M$ has five linear columns and five quadric columns, i.e.
$$M=[\quad L\quad | \quad Q \quad]$$ 
where $L$ is a $6\times 5$ matrix with linear forms and $Q$ is a $6\times5$ matrix with quadric forms.
 Also, by multiplying $M$ on left by the invertible matrix
\begin{equation}\label{eq:U}
U=
\left[ \begin{array}{cccccc}
1 & 0 & 0 & 0 & 0 & 0\\
0 & 1 & 0 & 0 & 0 & 0\\
0 & 0 & 1 & 0 & 0 & 0\\
0 & 0 & 0 & 1 & 0 & 0\\
0 & 0 & 0 & 0 & 1 & 0\\

c_1 & c_2 & c_3 & c_4 & c_5 & c_6
\end{array} \right]
\end{equation}
we make the last row of $L$ equal to $0$.
\end{proof}

\begin{proposition}
Algorithm \ref{Decomposition of cubic} computes the unique decomposition of general cubic $f$.
\end{proposition}
\begin{proof}
By Lemma \ref{correctness}, steps 1 through 3 of the Algorithm \ref{Decomposition of cubic} are well defined. For the step 4, consider the inverse of matrix $U$ defined in (\ref{eq:U}):
\begin{equation*}
V=\left[ \begin{array}{cccccc}
1 & 0 & 0 & 0 & 0 & 0\\
0 & 1 & 0 & 0 & 0 & 0\\
0 & 0 & 1 & 0 & 0 & 0\\
0 & 0 & 0 & 1 & 0 & 0\\
0 & 0 & 0 & 0 & 1 & 0\\

-c_1/c_6 & -c_2/c_6 & -c_3/c_6 & -c_4/c_6 & -c_5/c_6 & 1/c_6
\end{array} \right].
\end{equation*}
By applying the column operation given by the right multiplication of this matrix to the generators $[ \, g_1 \ldots g_6\,]$, we get the new generators 
\begin{equation*} 
[\,c_6 g_1 - c_1 g_6 \,c_6 g_2 - c_2 g_6 \ldots \,c_6 g_5 - c_5 g_6, g_6\,].
\end{equation*}
The first five of them have five linear syzygies given by the five columns of $UL$.  Moreover, choosing right basis as in \cite{MR} and applying Lemma \ref{correctness}, the second syzygy matrix $N$ also has the form 
$$N^t=[\quad Q'\quad |\quad L'\quad]$$
and the first row of $L'$ is $0$.

From this, the ideal $J$ generated by the first five generators $\,c_6 g_1 - c_1 g_6$, 
$\,c_6 g_2 - c_2 g_6\,$,$\ldots$, $\,c_6 g_5 - c_5 g_6$ has the minimal free resolution of the form
\begin{equation}\label{resolution of J}
0\longrightarrow S\longrightarrow S^5 \stackrel{P}{\longrightarrow} S^5 \longrightarrow S\longrightarrow 0
\end{equation}
where $P$ is given by the $5\times 5$ upper submatrix of $L$ as we seen above.

Since the resolution has the form (\ref{resolution of J}), the ideal $J$ has dimension $0$ and degree $5$, and it is also contained in the apolar ideal $f^\perp$. Further, by Lemma \ref{correctness}, the five chosen quadrics span the unique subspace contained in $f^\perp$ with this resolution.

Since there is the unique decomposition (\ref{eq:ffive}) for the given general cubic $f$ by \cite{Clebsch}, the linear forms in the decomposition define an ideal $K\subset f^\perp$ of five points in $\PP^3$. Since $K$ is Gorenstein of codimension 3, by the Buchsbaum-Eisenbud structure theorem, it has the free resolution of the form
\begin{equation}\label{five point resolution}
0\longrightarrow S\longrightarrow S^5 \stackrel{Q}{\longrightarrow} S^5 \longrightarrow S\longrightarrow 0.
\end{equation}
It means that the ideal $J$ we obtained by the Algorithm \ref{Decomposition of cubic} is actually same as the ideal $K$, which defines the five linear forms in (\ref{eq:ffive}), as we wanted.
\end{proof}
Using this algorithm, we can easily check whether the given quarternary cubic has real rank 5 or not. Namely, one computes the unique decomposition (\ref{eq:ffive}) and checks whether it is real. The real rank boundary can be obtained as in the following proposition. This proposition confirms \cite[Conjecture 5.5]{MMSV} for quaternary cubics.
\begin{proposition}\label{prop:R43}
The real rank boundary $\partial_{alg}(\mathcal{R}_{4,3})$ equals the join $J(\sigma_3(v_3(\mathbb{P}^3)),\tau(v_3(\mathbb{P}^3)))$ of third secant of the third Veronese embedding of $\mathbb{P}^3$ and its tangential variety. 
It is the irreducible hypersurface of degree $40$ in the $\PP(S^3(\CC^4)))$ with parametric representation
\begin{equation}\label{parametrization of join}g=\ell_1^3+\ell_2^3+\ell_3^3+\ell_4^2\ell_5, \qquad where \quad \ell_1,\ldots,\ell_5\in\RR[x,y,z,w]_1.
\end{equation}

\end{proposition}
\begin{proof}
The parametrization defines a unirational variety $Y$ in $\PP^{19}$. The Jacobian of this parametrization is found to have corank $1$. This means
that $Y$ has codimension $1$ in $\PP^{19}$.
 Hence $Y$ is an irreducible hypersurface, defined by a unique (up to sign)
 irreducible homogeneous polynomial $\Phi$ in $20$ unknowns with rational coefficients.
 
 Let $g$ be a real cubic with the form (\ref{parametrization of join}) that is a general point in $Y$. As $\epsilon$ goes to $0$, the real cubics $(\ell_4+\epsilon\ell_5)^3-\ell_4^3$ and $(i\ell_4+\epsilon \ell_5)^3+(-i\ell_4+\epsilon \ell_5)^3$ converge to the cubic $\ell_4^2\ell_5$ in $\PP^{19}$. It means that any small neighborhood of $g$ in $\PP^{19}$ contains cubics of real rank $5$ and cubics of real rank $> 5$. This implies that $Y$ lies in the real rank boundary $\partial_{\rm alg}\mathcal{R}_{4,3}$. Since $Y$ is irreducible and codimension $1$, it follows that $\partial_{\rm alg}\mathcal{R}_{4,3}$ exists and has $Y$ as an irreducible component.
 
 Using the algorithm \ref{Decomposition of cubic}, we can exactly compute the degree of the hypersurface $\partial_{\rm alg}\mathcal{R}_{4,3}$ which is $40$. This is done as follows. First, fix the field $K=\QQ(t)$ with a new variable $t$.  We fix two cubic $f_1$ and $f_2$ in $\QQ[x_1,x_2,x_3,x_4]_3$, and run the algorithm \ref{Decomposition of cubic} for $f_1+tf_2\in K[x_1,x_2,x_3,x_4]$. Step 4 returns an ideal $J$ in $K[x_1,x_2,x_3,x_4]$ that defines $5$ points in $\PP^3$ over the algebraic closure of $K$. By eliminating each of two variables of $\{x_1,x_2,x_3,x_4\}$, we obtain six binary forms of degree $5$ such that their coeffcients are degree $15$ polynomials in $t$. The discriminant of each binary form is a polynomial in $\QQ[t]$ of degree $120=8*15$. The greatest common divisor of these discriminants is a polynomial $\Psi(t)$ of degree $40$. We can check that $\Psi(t)$ is irreducible in $\QQ[t]$.

By definition, $\Phi$ is a
homogeneous polynomial with integer coefficients, irreducible over $\QQ$, in the $20$ coefficients
of a general cubic $f$. Its specialization $\Phi(f_1 + t f_2)$ is
a non-constant polynomial in $\mathbb{Q}[t]$, of degree ${\rm deg}(X)$ in $t$.
That polynomial divides $\Psi(t)$ because $Y$ lies in the real rank boundary. Since the latter is also irreducible over $\QQ$, we conclude that
 $\Phi(f_1 + t f_2) = \gamma \cdot \Psi(t)$, where $\gamma$ is
 a nonzero rational number. Hence $\Phi$ has degree $40$.
We conclude that ${\rm deg}(Y) = 40$, and therefore
$ Y = \partial_{\rm alg}( \mathcal{R}_{4,3})$.
\end{proof}

\subsection{Quaternary Quartics}
\begin{proposition}\label{prop:R44}
One component of the boundary $\partial_{alg}\mathcal{R}_{4,4}$ equals the dual of a $24$ dimensional variety of
quartic symmetroids in $\PP^3$, that is, the surfaces whose defining
polynomial is the determinant of a symmetric $4\times 4$-matrix of linear forms.
\end{proposition}
\begin{proof}
By Hilbert's classification we know that for quaternary quartics the cone of positive forms $P_{4,4}$ \emph{strictly} contains the cone of sums of squares $\Sigma_{4,4}$. Hence, the cone $Q_{4,4}$ of sums of (arbitrary many -- c.f.~Remark \ref{rem:overrank}) $4$-th powers of linear forms, which is the convex hull of the Veronese and equals the dual $P_{4,4}^*$, is strictly contained in the cone $\Sigma_{4,4}^*$ of forms with psd catalecticant. 

The boundary of $Q_{4,4}$ has two components. One (that is not interesting for our purposes) is the boundary of $\Sigma_{4,4}^*$ given by the determinant of the catalecticant. The other one, which we denote by $B$, is the dual of the Zariski closure of the extremal rays of $P_{4,4}\setminus\Sigma_{4,4}$ -- the variety of quartic symmetroids in $\PP^3$ by \cite[Theorem 3]{BHORS}.

The forms in $\Sigma_{4,4}^*\setminus Q_{4,4}$ are obviously of real rank greater than $10$, as they are not sums of powers by definition, and any other presentation would contradict the signature of the catalecticant. Further a generic point of $B$ has real rank at most $10$ by \cite[Proposition 7]{BHORS}. Changing the linear forms of the Waring decomposition of quartics in $B$ we obtain a Zariski dense set of forms of real rank at most $10$, intersecting $B$ in a relatively open set, which proves the proposition. 
\end{proof}
The following lemma is well-known to experts. We present a sketch of a proof based on \cite[Lemma 4.18]{Semidefinite}.
\begin{lemma}\label{positive rank}
For all postive integer $k$, the following set
$$C_{k,2d}:=\{~\sum_{i=1}^k\lambda_i\ell_i^{2d}\in S^{2d}(V^*)~|~\lambda_i\in\RR_{\geq 0},~ \ell_i\in V^*~\}$$
is a closed cone. In other words, the cone of forms whose positive rank is less than or equal to $k$ is closed.
\end{lemma}
\begin{proof}
 Let $\Delta=\{(\lambda_1,\dots,\lambda_k)\in \RR_{\geq 0}^k ~|~ \sum_{i=1}^k \lambda_i=1\}$ be a simplex and let $O$ be the unite sphere in $V^*$. Consider the image $C$ of the map:
 $$\Delta\times O^k\ni(\lambda_1,\dots,\lambda_k,\ell_1,\dots,\ell_k)\rightarrow \sum_{i=1}^k \lambda_i\ell_i^{2d}\in S^{2d}(V^*).$$
 Clearly, $C$ is compact and $0\not \in C$. Thus, the cone over $C$ is closed and it coincides with $C_{k,2d}$.
\end{proof}

\begin{remark}\label{rem:overrank}
There exists an nonempty open set inside $Q_{4,4}$ whose elements have real rank strictly greater than $10$. Indeed, the quartic $(x_1^2+x_2^2+x_3^2+x_4^2)^2$ has real rank $11$ by \cite[Proposition 9.26]{Rez} and hence belongs to $Q_{4,4}\setminus C_{10,4}$. The latter is open by Lemma \ref{positive rank}. Further, every element of this set has real rank strictly greater than $10$. Indeed, since every element of $Q_{4,4}$ has definite catalecticant matrix, it cannot have decompositions whose coefficients have distinct signs. 
\end{remark}   
One can easily prove that quaternary quartics $f$ of signature $(9,1)$ have rank at least $11$ if $\Omega(f)$ has no real points. This and Remark \ref{rem:overrank} motivate the following conjecture.
\begin{conjecture}
The real rank boundary $\partial_{alg}\mathcal{R}_{4,4}$ is reducible. The discriminant of $\Omega(f)$ is one of its components. Further, at least one more component comes from (the components of) the algebraic boundary of $C_{10,4}$. 
\end{conjecture}

\section{Appendix}\label{sec:App}
\subsection{Ternary Quadrics}\label{sec:AppTerQuad}
\begin{verbatim}
loadPackage("VarietyOfPolarSimplices")
p=0, n=3, (R,A,I)=unfoldingEquations(p,n)
J=flatteningRelations(R,A,I), S=ring I, X=(entries vars(R))_0
q=2*X_0*X_(n-1)+sum for i from 1 to n-2 list X_i^2
F=sub(I,S)+sub(J,S), B=sub(q,S), FintB=sub(F+B,sub(X_0,S)=>1)
P1=eliminate({sub(X_1,S),sub(X_2,S)}, FintB)
for i from 0 to 2 do
P2=sub(P1,sub(temp=leadMonomial(P1_i),S)=>sub(temp-1/leadCoefficient(P1_i)*P1_i,S));
Boundary=(gens P2)_3
\end{verbatim}
\subsection{Algebraic boundary $\partial_{alg}{\rm SSP}(xt+y^2+z^2)$}\label{app:Quat Quad}
\begin{verbatim}
loadPackage("VarietyOfPolarSimplices") 
p=0, n=4, (R,A,I)=unfoldingEquations(p,n)
J=flatteningRelations(R,A,I), S=ring I, X=(entries vars(R))_0
q=2*X_0*X_(n-1)+sum for i from 1 to n-2 list X_i^2
F=sub(I,S)+sub(J,S), B=sub(q,S), FintB=sub(F+B,sub(X_0,S)=>1)

K1=sub(FintB,sub(X_(n-1),S)=>sub(-1/2*(sum for i from 1 to n-2 list X_i^2),S))
KK1=K1
for i from 6 to 17 do
K1=sub(K1,sub(temp=leadMonomial(K1_i),S)=>sub(temp-1/leadCoefficient(K1_i)*K1_i,S))

K2=ideal(for i from 0 to 5 list K1_i)
Boundary=eliminate(K2,{sub(X_1,S),sub(X_2,S)})
\end{verbatim}
\subsubsection{$\partial_{alg}{\rm SSP}(xt+y^2+z^2)$ is prime.}\label{App:Boundary prime}
\begin{verbatim}
va=(entries sub(vars(A),S))_0
aa=for i from 6 to 17 list leadMonomial(KK1_i)
AA=QQ[toList(set(va)-set(aa)),t]
nB=homogenize(sub(Boundary,AA),t)
cB=ideal(random(1,AA), random(1,AA),random(1,AA),random(1,AA))+nB
eB=eliminate(cB,{(vars AA)_0_0,(vars AA)_1_0,(vars AA)_2_0,(vars AA)_3_0})
MM=QQ[a,b,t,c_0..c_40]
I=ideal(a,b,t)
gC=sub((gens gb eB)_0_0,matrix{{0,0,0,0,a,b,t}})
for i from 1 to 4 do
(
vars1=sub((i+2)*(i+1)/2,ZZ);
vars2=sub((8-i+2)*(8-i+1)/2,ZZ);
Am1=matrix{{c_0..c_(vars1-1)}};
Am2=matrix{{c_(vars1)..c_(vars1+vars2-1)}};
ff1=Am1*(transpose gens I^(i));
ff2=Am2*(transpose gens I^(8-i));
zer=sub(ff1*ff2-gC,{c_0=>1});
JJ=ideal(diff(gens I^8 , zer));
print degree JJ )
\end{verbatim}
\subsubsection{Real Points of $\mathfrak{Y}_4$ and $\mathfrak{Y}_{4,sing}$}\label{App:RealPoints}
\begin{verbatim}
sl=radical ideal singularLocus Boundary
y31=(primaryDecomposition(sl))_0, y22=(primaryDecomposition(sl))_1
y4=radical ideal singularLocus y31, y4sing=radical ideal singularLocus y22
f1=va_12-3*va_16,f2=va_17-3*va_14
y4_0==f1^2+f2^2
y4+ideal(f1,f2)==y4sing
\end{verbatim}

\subsection{Quinary Quadric}\label{App:QuinaryQuadrics}
\begin{verbatim}
loadPackage("VarietyOfPolarSimplices")
p=1009, n=5, (R,A,I)=unfoldingEquations(p,n)
J=flatteningRelations(R,A,I), S=ring I
X=(entries vars(R))_0
q=2*X_0*X_(n-1)+sum for i from 1 to n-2 list X_i^2
F=sub(I,S)+sub(J,S), B=sub(q,S)
FintB=sub(F+B,sub(X_0,S)=>1);
K2=sub(FintB,sub(X_(n-1),S)=>sub(-1/2*(sum for i from 1 to n-2 list X_i^2),S));
K3=K2;
for i from 10 to 39 do
K3=sub(K3,sub(temp=leadMonomial(K3_i),S)=>sub(temp-1/leadCoefficient(K3_i)*K3_i,S))
K4=ideal(for i from 0 to 9 list K3_i);
va=(entries sub(vars(A),S))_0, aa=for i from 10 to 39 list leadMonomial(K2_i)
SS=ZZ/1009[toList(set(va)-set(aa)),X,t]
K4=sub(K4,SS);
K4=ideal(homogenize(gens K4,t));
for i from 0 to 7 do 
K4=sub(K4,{SS_i=>random(ZZ/1009)*SS_8+random(ZZ/1009)*SS_9+random(ZZ/1009)*t});
L=eliminate(K4,{sub(X_1,SS),sub(X_2,SS),sub(X_3,SS)})
L=saturate(L,t)
degree L
\end{verbatim}

\subsection{Quarternary Cubic}\label{sec:AppQuartCub}
\begin{verbatim}
A=QQ[a_0..a_3]
B=A[b_(0,0)..b_(5,3)]
R=B[x_1,x_2,x_3,x_4]
bs:=d->(entries basis(d,R))_0;
f=a_0*x_1^3+a_1*x_1^2*x_2+a_2*x_1^2*x_3+a_3*x_1^2*x_4+11*x_1*x_2^2-12*x_1*x_2*x_3+
7*x_1*x_2*x_4+32*x_1*x_3^2-28*x_1*x_3*x_4+11*x_1*x_4^2+8*x_2^3-13*x_2^2*x_3+
34*x_2^2*x_4+19*x_2*x_3^2-38*x_2*x_3*x_4+16*x_2*x_4^2+7*x_3^3-41*x_3^2*x_4+
7*x_3*x_4^2+13*x_4^3;
cat=sub(diff(transpose basis(1,R)*basis(2,R),f),A)
kn=kernel cat
gsI=for i from 0 to 5 list sum for j from 0 to #(bs(2))-1 list kn_i_j*(bs(2))_j
I=ideal gsI
candSyz=for i from 0 to 5 list sum for j from 0 to 3 list b_(i,j)*R_j
allsum=sum for i from 0 to 5 list candSyz_i*gsI_i
bb=flatten for i from 0 to 5 list for j from 0 to 3 list b_(i,j)
MM=transpose matrix for j from 0 to #bb-1 list 
for i from 0 to #(bs(3))-1 list coefficient(bb_j,coefficient((bs(3))_i,allsum))
rank MM
for i from 1 to 4 list rank submatrix(MM,0..19,i..19+i)
rank submatrix(MM,0..19,0..0|2..20)
kn1to4=for i from 1 to 4 list(
    temp=kernel submatrix(MM,0..19,i..19+i);
    transpose matrix{for j from 0 to 23 list 
    if(j<i) then 0 else if(i<=j and j<20+i) then temp_0_(j-i) else 0}
    );
knMatrix=kn1to4_0;
for i from 1 to 3 do knMatrix=knMatrix|kn1to4_i;
temp=kernel (sm=submatrix(MM,0..19,0..0|2..20));
kn5=transpose matrix{for j from 0 to 23 list 
if(j==0) then temp_0_0 else if(j==1) then 0 
else if(j>1 and j<=20) then temp_0_(j-1) else 0};
knMatrix=knMatrix|kn5;
rank knMatrix
linSyz=for k from 0 to 4 list for i from 0 to 5 list 
sum for j from 0 to 3 list knMatrix_k_(i*4+j)*R_j;
L=transpose matrix{linSyz_0};
for i from 1 to 4 do L=L|transpose matrix{linSyz_i};

C=A[c_0..c_5]
S=C[x_1,..,x_4]
LC=sub(L,S);
linSum=for j from 0 to 4 list sum for i from 0 to 5 list c_i*LC_j_i;
knlist=for k from 0 to #linSum-1 list kernel transpose matrix for j from 0 to 5 list 
for i from 0 to 3 list coefficient(c_j,coefficient(S_i,linSum_k));
clist=intersect knlist;
rank clist
J=ideal(for i from 0 to 4 list gsI_i*clist_0_(i+1)-gsI_(i+1)*clist_0_i);
\end{verbatim}
\medskip

\begin{small}

\end{small}

\bigskip

\noindent
\footnotesize {\bf Authors' addresses:}

\smallskip

\noindent Mateusz Micha{\l}ek,
Institute of Mathematics, 
Polish Academy of Sciences, Warsaw, Poland,\\
Max Planck Institute MiS, Leipzig
{\tt wajcha2@poczta.onet.pl}

\smallskip

\noindent Hyunsuk Moon,
  KAIST,
     Daejeon, South Korea, {\tt octopus14@kaist.ac.kr}

\end{document}